\documentclass{amsart}
\usepackage[foot]{amsaddr}
\usepackage{amsmath}
\usepackage{amssymb} 
\usepackage{graphicx} 

\numberwithin{equation}{section}

\def\ad#1{\begin{aligned}#1\end{aligned}}  \def\b#1{\mathbf{#1}} 
\def\a#1{\begin{align*}#1\end{align*}} \def\an#1{\begin{align}#1\end{align}} \def\t#1{\hbox{#1}}

\begin{document} 

\title[finite elements]{ A nodal basis for the $C^1$-$P_{33}$ 
    finite elements on 5D simplex grids
}

 \author {Jun Hu}
\address{LMAM and School of Mathematical Sciences, Peking University, Beijing 100871, P. R. China.}
\email{hujun@math.pku.edu.cn}

 \author {Shangyou Zhang}
\address{Department of Mathematical Sciences, University of Delaware,
    Newark, DE 19716, USA. }
\email{szhang@udel.edu}

\begin{abstract}
We construct a nodal basis for the 5-dimensional $C^1$ finite element space of polynomial degree $33$
  on simplex grids,
  where the finite element functions are $C^1$ on the 6 4D-simplex faces,
     $C^2$ on the 15 face-tetrahedra, $C^4$ on the 20 face-triangles, 
     $C^8$ on the 15 edges, and $C^{16}$ at the 6 vertices, of a 5D simplex.
     
\end{abstract}



\maketitle 

 \baselineskip=16pt

\section{Introduction} 
In 2D, the Argyris finite element is one of the first finite elements even constructed,
     cf. \cite{Ciarlet}.
It is a $C^1$-$P_5$ finite element on triangular grids.  
Here $C^1$-$P_5$ denotes the space of globally $C^1$ and locally piecewise polynomials of
   degree 5 on 2 dimensional triangular grids. 
But the Argyris element is $C^2$ at vertices and $C^1$ at edges, with the degrees of freedom 
   of up-to order 2 derivatives at each vertex and 1 normal derivative at each edge.
It is straightforward to extend the Argyris finite element to $C^1$-$P_{k}$ ($k> 5$) finite elements,
    as follows.
We introduce $(k-5)$ function values at $(k-5)$ internal points on each edge; we 
    introduce first-order normal derivatives at
            $(k-4)$ internal points each edge;
  and we introduce 
           function values at 2D $\dim P_{k-6}$ internal points on the triangle.

Also in 2D, in 1970, Bramble and Zl\'{a}mal \cite{Bramble}   and  \v{Z}en\'i\v{s}ek \cite{Zenisek}
  extended the above $C^1$-$P_{k}$ finite element to $C^m$-$P_{4m+1}$ finite elements
  for all $m\ge 1$, on triangular grids.
In fact, \v{Z}en\'i\v{s}ek \cite{Zenisek} defined 
   all $C^m$-$P_{k}$ finite elements for $k\ge 4m+1$.
   
In 3D, the first $C^1$ element was constructed
     by  \v{Z}en\'i\v{s}ek in 1973 \cite{Zenisek-3d}, a $C^1$-$P_9$
   finite element.
By avoiding high-order derivatives in the degrees of freedom of \v{Z}en\'i\v{s}ek  \cite{Zenisek-3d},
   the element was extended  to all $C^1$-$P_k$  ($k\ge 9$) finite elements \cite{Zhang-3d}
   in 2009.
Also \v{Z}en\'i\v{s}ek  extended the $C^1$-$P_9$ finite element to $C^m$-$P_{8m+1}$ on tetrahedral grids
   in 1974 \cite{Zenisek-cm}.
Again high-order derivatives (above the continuity order) were used as degrees of freedom in \cite{Zenisek-cm}.
That is, the \v{Z}en\'i\v{s}ek 3D $C^1$-$P_9$ and 3D $C^2$-$P_{17}$ are not super-spline finite elements.
A different $C^2$-$P_{k}$ ($k\ge 17$) on tetrahedral grids was constructed in \cite{Zhang-4d}, in 2016.
The family of $C^m$-$P_{k}$ ($k\ge 8m+1$) for any $m\ge 1$
   on tetrahedral grids was constructed in \cite{Zhang-34d}, in 2022.

In 4D, a family of $C^1$-$P_{k}$ ($k\ge 17$) finite elements on 4D-simplex grids was
  constructed in \cite{Zhang-4d}, in 2016.
The element was extended to $C^m$-$P_{k}$ ($k\ge 16m+1$) for any $m\ge 1$
   on 4D simplex grids in \cite{Zhang-34d}, in 2022.
 
In 5D, we construct, in this work, a nodal basis for $C^1$-$P_{33}$ finite element on 5D simplex grids,
  which is $C^{16}$ at vertices, $C^{8}$ at edges, $C^4$ at triangles, $C^2$ at tetrahedra and
    $C^1$ on face 4D-simplexes.

Such finite elements are called super-splines for excessive smoothness at low-dimensional simplex.
Alfeld, Schumaker and Sirvent are the first to
   introduce a distance function from a Bernstein index
     to a low-dimensional simplex and to define the Bernstein index partition for the nodal basis,  
      cf.  Equation (36) in \cite{Alfeld}.
As claimed by \cite{Alfeld}, it is very difficult to find an 
   explicit partition of these index sets for general or 
  a given set of $m$, $n$ and $k$ (for the smoothness, the space dimension and the polynomial degree). 
In addition to \cite{Alfeld}, \cite{Chen} and \cite{Hu}  proposed and proved the index partition.

Mathematically it is a challenge to find explicit definitions of basis functions for 
    general $C^m$-$P_{k}$ finite elements in high-space dimensions.
As we can see, in this work, the Bernstein index partition has more and more overlapping
  in higher and higher space dimensions.
Only in 2D, there is no overlapping for all $C^m$ finite elements.

\section{The vertex, edge and triangle degrees of freedom}

Let $\mathcal{T}_h=\{K\}$ be a 5D simplex grid where a simplex 
   $K$ has 6 vertices $\{\b x_i=
   (x_1^{(i)}, x_2^{(i)}, x_3^{(i)}, x_4^{(i)}, x_5^{(i)}), \ i=1,\dots,6\}$.
The intersection of
   of two $K$'s is either a  4D common face-simplex, or a common face-tetrahedron,
   or a common triangle, or a common edge, or a common vertex, or an empty set.
The maximal diameter of $\{K\}$ in $\mathcal T_h$ is $h$.

\long\def\mskip#1{}

The vertex degrees of freedom have the highest priority that it is required all derivatives
   from order 0 to 16 to be specified at each vertex.
By the $C^{16}$ requirement, we define the following degrees of freedom at the six vertices of
  a 5D simplex,
\an{\label{b11}  \partial_{\b x^{\alpha} } p(\b x_i), \quad |\alpha|\le 16, \ i=1,\dots,6, }
where $\alpha=\langle \alpha_1,\dots, \alpha_5\rangle$ is from the Bernstein index
\a{  (\alpha_1,\dots, \alpha_5, 33 - |\alpha|), \quad |\alpha|=  \alpha_1+\cdots+\alpha_5\ge 0.  }
For example, the Bernstein index is $(0,0,0,1,15,17)$ for the 16-th derivative nodal basis function
\a{ \partial_{\b x_4\b x_5^{15} } \phi_j(\b x_i)=\delta_{ij}. } 

Summing up the number of derivatives (degrees of freedom) in \eqref{b11}, we have
\a{ &\quad \ C_{5,0}+C_{5,1}+C_{6,2}+\dots + C_{19,15}+ C_{20,16}\\
          &= 1 + 5 + 15 + \dots + 3876+  4845  \\
          &= C_{21,16}=\frac{ 21! }{16! \cdot 5! }=20349, }
where 
\a{ C_{m,n} = \frac{ m ! } { (m-n)! \cdot n! }, }
and the following hockey-stick identity is applied, 
\a{ \sum_{i=0} ^ k C_{n+i,i} = C_{n+k+1,k}.  }
But we will use another form of this hockey-stick identity repeatedly, to combine
  the counting on low-dimensional simplexes to one on one-dimensional higher simplex,
   that
\a{ \sum_{i=0} ^ k C_{n+i,n} = C_{n+k+1,n+1}.  }

As there are 6 vertices on a 5D simplex $K$,  the total number of vertex degrees of freedom is
\an{\label{1-s} \t{dof}_{\t{vert}}  &= 6\cdot 20349 =  122094.  }

For a $P_{33}$ polynomial, after specifying its 2 end-point derivatives up to order 16,
  we have $2(16+1)=34$ constraints on its 34 coefficients on the 1D polynomial from the restriction
    of the 5D polynomial on the edge.
Thus, we have 0 point value, 4 first normal derivatives at 1 point on the edge
   (with the Bernstein indexes $(1,0,0,0,16,16)$--$(0,0,0,1,16,16)$), 
   10 second normal derivatives at 2 points, and so on until 165 8-th normal derivatives at
   8 points, 
\an{\label{b-2} \partial_{\b x^{\alpha}} (\b m_{j}), \ 0 <|\alpha|\le 8, \ \alpha_5=0,
     \ j=0,\dots,8, \  }
where $\b x$ is understood orthogonal to one $\b x_i$, and
   $\b m_j$ are 1D $P_j$ Lagrange nodes inside an edge.
Here we drop $\b x_i$ to denote an orthogonal direction to $\b x_i\b x_6$, which means, 
 in other words, $\b x_6=\langle 0,0,0,0,0\rangle$ and $\b x_1=\langle 1,0,0,0,0\rangle$.
 
 We sum all edge degrees of freedom in \eqref{b-2} to get 
\a{ & \quad \  0+  1\cdot C_{4,3}+2 C_{5,3} + 3 C_{6,3} + \dots + 8 C_{8,3} \\
    & =  0+4+2\cdot 10+3\cdot 20 +\cdot + 8\cdot 165 \\
    &= 3168. }
As there are 15 edges on a $K$,  the total number of vertex degrees of freedom is
\an{\label{2-s}   \t{dof}_{\t{edge}} &= 15 \cdot 3168 =  47520.  }

\mskip{ 
restart; ft[0]:=1: for i to 25 do ft[i]:=i*ft[i-1]: od: # for i from 0 to 20 do i,ft[i]; od;
 for i from 1 to 16 do i, ft[5+i-1]/(ft[i]*ft[4]),ft[5+i]/(ft[5]*ft[i]); od;

for i to 8 do i, i*ft[3+i]/(ft[i]*ft[3]); od;

for i from 0 to 4 do i,[3+i,i], ft[3+i]/(ft[i]*ft[3]); od;
 
i:=6: j:=1: ft[i]/(ft[i-j]*ft[j]);

}

For the face-triangle degrees of freedom, to be $C^4$ continuous,  we need additional 
  2D-$P_6$ function values, 2D-$P_9$ first normal derivatives in each of three ($C_{3,2}$) directions,
  2D-$P_{12}$ second normal derivatives in  $C_{4,2}$ directions,
  2D-$P_{15}$ third normal derivatives in  $C_{5,2}$ directions,
  and  2D-$P_{18}$ 4-th normal derivatives in $C_{6,2}=15$ directions.
  
As we defined up-to $P_{16}$ degrees of freedom of zero at two end-points and up-to $P_8$ 
  degrees of freedom of zero at the edge,  we can factor $\lambda_1^9$ from the function where
   $\lambda_1$ is the bary-centric coordinate of the edge on a triangle.
Thus the function $p$ is of the form $p=(\lambda_1\lambda_2\lambda_3)^9 p_{6}$, where
  $p_6$ is a $P_6$ polynomial on the face-triangle.
We have a perfect match here, only for function values on 2D $P_6$ internal points.
We will have overlapping of index for the first normal derivative and so on.
To determine the $P_6$ polynomial, and also to make the function $C^0$ continuous,  we 
  need the following $C_{6,2}=28$ degrees of freedom,
\a{ p(\b x_{\ell,i}^{(2)}), \quad \ell=1,\dots, C_{6,3}=20, \ i=1,\dots, C_{8,2}=28, }
where $\b x_{\ell,i}^{(2)}$ are internal Lagrange $P_6$ points on one of $20$ face-triangles.
The number of function-value degrees of freedom inside all triangles is
\an{\label{triangle-0} \t{dof}_{\t{tri},0} = 20\cdot 28 = 560. }

Typically, we need to specify the first normal derivatives at
     2D $P_{6+3}$ internal Lagrange points.
However, for a normal derivative, we have a $P_{32}$ polynomial for the $P_{33}$ function $p$.
After specifying its up-to order $16-1$ derivatives at the three corners and its up-to order $8-1$ 
  normal derivatives at the three edges,  we have only a set of 2D
    $P_{32-3\times 8}=P_8$ Lagrange nodes left
    for the wanted 2D $P_9$ polynomial space above.
That is, the face-normal-derivative required indexes $(1,0,0,7,*,*)$ are already
  covered by the vertex and edge degree of freedom above.
Therefore, we do not need to specify these 10 ($=\dim P_9-\dim P_8$ in 2D) normal derivatives
  as the total $\dim P_{32}$(in 2D) $C^1$-constraints are already specified on the face triangle.
 That is,  we need the following $C_{10,2}=45$($=\dim P_8$ in 2D) degrees of freedom,
\a{ \partial{\b x_{(3)} } p(\b x_{\ell,i}^{(2)}), \quad \ell=1,\dots, C_{6,3}=20, \ i=1,\dots, C_{10,2}=45, }
where $\partial{\b x_{(3)}}$ denotes the first derivatives in three normal directions,
  and $\b x_{\ell,i}^{(2)}$ are internal Lagrange 2D-$P_8$ points on one of $20$ face-triangles.
The number of first-normal-derivative degrees of freedom inside all triangles is
\an{\label{triangle-1} \t{dof}_{\t{tri},1} =3\cdot 20\cdot 45 = 2700. }

For the second derivative, we have a 2D $P_{31}$ polynomial on a face-triangle.
As the three edges' restricted by 1D $P_6$ polynomials, we have
  a 2D $P_{31-3\times 7}=P_{10}$ polynomial degrees of freedom to be specified,
   not the 2D $P_{12}$  polynomial degrees of freedom mentioned above.
However, unlike the first normal derivatives,
    as we had $P_{16-2}$ derivatives' restrictions at three vertices,
  the Bernstein index $(2,0,0,7,7,17)$ for a corner $P_{10}$ is already covered by the
    vertex degrees of freedom.
Thus,   we need to specify the second normal derivatives in
   2D $\dim P_{10}-3=66-3=63$ internal Lagrange points.
That is,  we need the following $C_{12,2}-3=63$ degrees of freedom for each of $C_{4,2}=6$ second
  normal derivatives,
\a{ \partial{\b x_{(3)}^2 } p(\b x_{\ell,i}^{(2)}), \quad \ell=1,\dots, 20, \ i=1,\dots, 63, }
where $\partial{\b x_{(3)}^2}$ denotes the order two derivatives in six normal directions,
  and $\b x_{\ell,i}^{(2)}$ are internal Lagrange 2D-$P_10$ minus 3 times 2D-$P_0$
      points at three vertices, on one of $20$ face-triangles.
The number of second-normal-derivative degrees of freedom inside all triangles is
\an{\label{triangle-2} \t{dof}_{\t{tri},2} =6\cdot 20\cdot 63 = 7560. }

There are $C_{2+3,2}=10$ third order directional normal derivatives on a face-triangle.
Again, due to overlapping, we 
  need to post $P_{12}$ (instead of $P_{15}$ mentioned above)
  2D Lagrange nodal degrees of freedom to the third derivative.
But the Bernstein index $(3,0,0,6,6,18)$, index $(3,0,0,6,7,17)$ and index $(3,0,0,7,6,17)$, are
  covered by the degrees of freedom at the vertex $\b x_6$. 
That is,  we need the following $C_{14,2}-3\times 3=82$ degrees of freedom for each of $C_{5,2}=10$
  third normal derivatives,
\a{ \partial{\b x_{(3)}^3 } p(\b x_{\ell,i}^{(2)}), \quad \ell=1,\dots, 20, \ i=1,\dots, 82, }
where $\partial{\b x_{(3)}^3}$ denotes the order three derivatives in ten normal directions,
  and $\b x_{\ell,i}^{(2)}$ are internal Lagrange 2D-$P_{12}$ minus 3 times 2D-$P_1$
      points at three vertices, on one of $20$ face-triangles.
The number of third-normal-derivative degrees of freedom inside all triangles is
\an{\label{triangle-3} \t{dof}_{\t{tri},3} =10\cdot 20\cdot 82  = 16400. }

There are $C_{2+4,2}=15$ fourth order directional normal derivatives on a face-triangle.
Supposedly we need to post $\dim P_{14}=120$ (instead of $P_{18}$ mentioned above)
   2D Lagrange nodal degrees of freedom to the fourth derivative
  of the function in order to have a $C^4$ continuity crossing the face-triangle.
But the Bernstein index $(4,0,0,5,5,19)$, index $(4,0,0,5,6,18)$ and index $(4,0,0,5,7,17)$, are
  covered by the degrees of freedom at the vertex $\b x_6$. 
That is,  we need only the following $C_{16,2}-3\times C_{4,2}=102$ 
  degrees of freedom for each of $C_{6,2}=15$ fourth normal derivatives,
\a{ \partial{\b x_{(3)}^4 } p(\b x_{\ell,i}^{(2)}), \quad \ell=1,\dots, 20, \ i=1,\dots, 102, }
where $\partial{\b x_{(3)}^4}$ denotes the order three derivatives in ten normal directions,
  and $\b x_{\ell,i}^{(2)}$ are internal Lagrange 2D-$P_{14}$ minus 3 times 2D-$P_2$
      points at three vertices, on one of $20$ face-triangles.
The number of fourth-normal-derivative degrees of freedom inside all triangles is
\an{\label{triangle-4} \t{dof}_{\t{tri},4} =15\cdot 20\cdot 102 = 30600. }

Together, by \eqref{triangle-0}, \eqref{triangle-1}, \eqref{triangle-2}, \eqref{triangle-3} and
   \eqref{triangle-4}, we have the following number of degrees of freedom inside all triangles, 
\an{\label{tri} \ad{ \t{dof}_{\t{tri}} &=
 \t{dof}_{\t{tri},0} + \t{dof}_{\t{tri},1}  + \t{dof}_{\t{tri},2}  + \t{dof}_{\t{tri},3}
    + \t{dof}_{\t{tri},4}  \\
  &= 560+2700+7560+ 16400+30600=57820.  } }

\mskip{ 
 check:    9   9  15   0   0       0   3   1   1   1   2   3
check:    8   8  16   0   0       1   3   2   1   1   2   3
check:    7   8  16   0   0       2   3   3   1   1   2   3
check:    6   8  16   0   0       3   3   4   1   1   2   3
check:    5   8  16   0   0       4   3   5   1   1   2   3
simplex 0  derivative 0 dof       1  sum=       1
simplex 0  derivative 1 dof       5  sum=       6
simplex 0  derivative 2 dof      15  sum=      21
simplex 0  derivative 3 dof      35  sum=      56
simplex 0  derivative 4 dof      70  sum=     126
simplex 0  derivative 5 dof     126  sum=     252
simplex 0  derivative 6 dof     210  sum=     462
simplex 0  derivative 7 dof     330  sum=     792
simplex 0  derivative 8 dof     495  sum=    1287
simplex 0  derivative 9 dof     715  sum=    2002
simplex 0  derivative10 dof    1001  sum=    3003
simplex 0  derivative11 dof    1365  sum=    4368
simplex 0  derivative12 dof    1820  sum=    6188
simplex 0  derivative13 dof    2380  sum=    8568
simplex 0  derivative14 dof    3060  sum=   11628
simplex 0  derivative15 dof    3876  sum=   15504
simplex 0  derivative16 dof    4845  sum=   20349
level   0  #simplex   6 dofs  20349 total  122094
simplex 1  derivative 0 dof       0  sum=       0
simplex 1  derivative 1 dof       4  sum=       4
simplex 1  derivative 2 dof      20  sum=      24
simplex 1  derivative 3 dof      60  sum=      84
simplex 1  derivative 4 dof     140  sum=     224
simplex 1  derivative 5 dof     280  sum=     504
simplex 1  derivative 6 dof     504  sum=    1008
simplex 1  derivative 7 dof     840  sum=    1848
simplex 1  derivative 8 dof    1320  sum=    3168
level   1  #simplex  15 dofs   3168 total   47520
simplex 2  derivative 0 dof      28  sum=      28
simplex 2  derivative 1 dof     135  sum=     163
simplex 2  derivative 2 dof     378  sum=     541
simplex 2  derivative 3 dof     820  sum=    1361
simplex 2  derivative 4 dof    1530  sum=    2891
level   2  #simplex  20 dofs   2891 total   57820

simplex 3  derivative 0 dof     544  sum=     544
simplex 3  derivative 1 dof    1778  sum=    2322
simplex 3  derivative 2 dof    3804  sum=    6126

level   3  #simplex  15 dofs   6126 total   91890
simplex 4  derivative 0 dof    6965  sum=    6965
simplex 4  derivative 1 dof   12990  sum=   19955
level   4  #simplex   6 dofs  19955 total  119730
simplex 5  derivative 0 dof   62888  sum=   62888
level   5  #simplex   1 dofs  62888 total   62888
 (n m k_1)= 5 1 0, dim P_{ 33}=  501942 C^m-P_k^n=  501942
}
 
 \section{The tetrahedral degrees of freedom }

 As the function is $C^4$-continuous on the 4 face-triangles of a tetrahedron,
supposedly we need to post $\dim P_{33-4\times 5}=\dim P_{13}=C_{16,3}=560$ 
   3D Lagrange nodal degrees of freedom to the function
   values in order to have a $C^0$ continuity crossing a face-tetrahedron.
But the Bernstein index $(0,0,5,5,5,18)$ and index $(0,0,5,5,6,17)$ are
  covered by the degrees of freedom at a vertex. 
That is,  we need only the following $C_{16,3}-4\times C_{6,3}=544$ 
  degrees of freedom for the function value on a face-tetrahedron,
\a{  p(\b x_{\ell,i}^{(3)}), \quad \ell=1,\dots, 15, \ i=1,\dots, 544, }
where   $\b x_{\ell,i}^{(3)}$ are internal Lagrange 3D-$P_{13}$ minus 4 times 3D-$P_3$
      points at four vertices, on one of $15$ face-tetrahedra.
The number of degrees of freedom for the function value inside all tetrahedra is
\an{\label{tet-0} \t{dof}_{\t{tet},0} =1\cdot 15\cdot 544 = 8160. }

On a face-tetrahedron,  we have 2 normal derivatives in 2 directions.
 As the first normal derivative of $p$
    is $C^3$-continuous on the 4 face-triangles of a tetrahedron,
supposedly we need to post $\dim P_{32-4\times 4}=\dim P_{16}=C_{19,3}=969$ 
   3D Lagrange nodal degrees of freedom to the first normal derivatives
      in order to have a $C^1$ continuity crossing a face-tetrahedron.
But the Bernstein index $(1,0,4,4,4,20)$, index $(1,0,4,4,5,19)$, index $(1,0,4,4,6,18)$ and
    index $(1,0,4,4,7,17)$ are
  covered by the degrees of freedom at a vertex. 
That is,  we need only the following $C_{19,3}-4\times C_{6,3}=889$ 
  degrees of freedom for the function value on a face-tetrahedron,
\a{ \partial{\b x_{(3)}^1 }  p(\b x_{\ell,i}^{(3)}), \quad \ell=1,\dots, 15, \ i=1,\dots, 889, }
where $\partial{\b x_{(3)}^1}$ denotes the order one derivatives in two normal directions,
   and   $\b x_{\ell,i}^{(3)}$ are internal Lagrange 3D-$P_{16}$ minus 4 times 3D-$P_3$
      points at four vertices, on one of $15$ face-tetrahedra.
The number of degrees of freedom for the normal derivative inside all tetrahedra is
\an{\label{tet-1} \t{dof}_{\t{tet},1} =2\cdot 15\cdot 889 = 26670. }

On a face-tetrahedron,  we have 3 second normal derivatives in 3 directions.
 As the second normal derivative of $p$
    is $C^2$-continuous on the 4 face-triangles of a tetrahedron,
supposedly we need to post $\dim P_{31-4\times 3}=\dim P_{19}$ 
   3D Lagrange nodal degrees of freedom for the second normal derivatives
      in order to have a $C^2$ continuity crossing a face-tetrahedron.
But the above 3D-$P_{19}$ constraints are based on the restrictions from vertices and edges.
In this case, restrictions from other face-triangles of the tetrahedron define the
  three outside lines of degrees of freedom on this face.
That is, the value of Bernstein index $(2,0,2,*,*,*)$ on this face is defined previously on 
   other face-triangles.
Hence, on the first level of the this face of the tetrahedron,  we have a 2D-$P_{16}$ polynomial
   instead of a 2D-$P_{19}$ polynomial.
As the Bernstein index $(2,0,3,4,7,17)$, index $(2,0,3,4,6,18)$, index $(2,0,3,4,5,19)$
  and index $(2,0,3,4,4,20)$ are covered by the vertex degrees of freedom, 
 we need only the following $C_{18,2}-3\times C_{5,2}=123$ 
  degrees of freedom for the second normal derivatives on the first level triangle of this
     face-tetrahedron,
\an{\label{t-0-0} \partial{\b x_{(3)}^2 }  p(\b x_{\ell,i}^{(3)}), \quad \ell=1,\dots, 1, \ i=1,\dots, 123, }
where $\partial{\b x_{(3)}^2}$ denotes the order two derivatives in three normal directions,
   and   $\b x_{\ell,i}^{(3)}$ are internal Lagrange 2D-$P_{16}$ minus 3 times 2D-$P_3$
      points at four vertices, on one face-triangle of this face-tetrahedra.
      
Since dofs on the tetrahedron has a pattern after chopping off 4-face triangles in \eqref{t-0-0}, 
    we will count the rest dofs in 2D any more.
   
After taking off the above degrees of freedom on the four face-triangles in \eqref{t-0-0},  we have
  a 3D-$P_{19-4}=P_{15}$ Lagrange nodes of degrees of freedom inside the tetrahedra.
But the Bernstein index $(2,0,4,4,7,17)$,   $(2,0,4,4,6,18)$,   $(2,0,4,4,5,19)$, 
  and   $(2,0,4,4,4,20)$ are still covered by the vertex degrees of freedom.
Thus, inside the tetrahedron, 
 we need additionally the following $C_{18,3}-4\times C_{6,3}=776$ 
  degrees of freedom for the second normal derivatives,
\a{ \partial{\b x_{(3)}^2 }  p(\b x_{\ell,i}^{(3)}), \quad \ell=1,\dots, 1, \ i=1,\dots, 776, }
where $\partial{\b x_{(3)}^2}$ denotes the order two derivatives in three normal directions,
   and   $\b x_{\ell,i}^{(3)}$ are internal Lagrange 3D-$P_{15}$ minus 4 times 3D-$P_3$
      points at four vertices.
Together, the
    number of degrees of freedom for the order-two normal derivative inside all tetrahedra is
\an{\label{tet-2} \t{dof}_{\t{tet},2} =3\cdot 15 ( 4\cdot 123 + 776 ) =  57060. }

Together, by \eqref{tet-0}, \eqref{tet-1} and \eqref{tet-2}, the total number of
   tetrahedral degrees of freedom inside the 15 face-tetrahedra is
\an{\label{tet} \ad{ \t{dof}_{\t{tet}} 
 &= \t{dof}_{\t{tet},0} + \t{dof}_{\t{tet},1} + \t{dof}_{\t{tet},2} \\
          &= 8160+26670 +57060 =91890.  } }

 \section{The 4D-simplex degrees of freedom }

 As the function is $C^2$-continuous on the 5 face-tetrahedra of a 4D-simplex,
  supposedly we need to post $\dim P_{33-5\times 3}=\dim P_{18}=C_{22,4}$ 
   4D Lagrange nodal degrees of freedom to the function
   values in order to have a $C^0$ continuity crossing the 4D-simplex.
But the Bernstein index $(0,3,3,3,7,17)$, index $(0,3,3,3,6,18)$,
    index $(0,3,3,3,5,19)$, index $(0,3,3,3,4,20)$ and index $(0,3,3,3,3,21)$ are
       covered by the degrees of freedom at a vertex. 
That is,  we need only the following $C_{22,4}-5\times C_{8,4}=6965$ 
  degrees of freedom for the function value on a face-tetrahedron,
\a{  p(\b x_{\ell,i}^{(4)}), \quad \ell=1,\dots, 6, \ i=1,\dots, 6965, }
where   $\b x_{\ell,i}^{(4)}$ are internal Lagrange 4D-$P_{18}$ minus 5 times 4D-$P_4$
      points at four vertices, on one of $6$ face-4-simplex.
The number of degrees of freedom for the function value inside all 4D-simplex is
\an{\label{s4-0} \t{dof}_{\t{4s},0} =1\cdot 6\cdot 6965 = 41790. }

On a 4D-simplex,  we have only one normal derivative in the normal direction.
 As the first normal derivative of $p$
    is $C^1$-continuous on the 5 face-tetrahedra of a 4D-simplex,
supposedly we need to post $\dim P_{32-5\times 2}=\dim P_{22}$ 
   4D Lagrange nodal degrees of freedom to the first normal derivatives
      in order to have a $C^1$ continuity crossing a face-tetrahedron.
We consider the degrees of freedom on the first layer of tetrahedra on the surface of this 4D simplex.
On the first triangle of this tetrahedron,  as restricted by $C^8$ edge-continuity,
  we have only 2D $P_{16}$ (instead of $P_{22}$ mentioned above)
    degrees of freedom left, at most,
  ranging from Bernstein index $(1,2,2,4,8,16)$ to index $(1,2,2,4,16,8)$ on the first edge of the
   triangle.
Thus, on this triangle, we have 2D 
\an{\label{2d-1} R^{(2)}_{16,-3} = \dim P_{16}-3 \dim P_3=123 }
  degrees of freedom as
  the index $(1,2,2,4,7,17)$, $(1,2,2,4,6,18)$,  $(1,2,2,4,5,19)$ and $(1,2,2,4,4,20)$ are
   covered by the vertex degrees of freedom. 
Here we introduce a short notation in \eqref{2d-1} for chopping corners that
\a{ R^{(k)}_{m,-n} = \dim P_m - (k+1)\dim P_n = C_{m+k,k}-(k+1)C_{n+k,k}.  }
Thus, from \eqref{2d-1},  we need only the following  
  degrees of freedom for the normal derivatives on the first triangle of the first tetrahedron
      of this 4D-simplex,
\a{ \partial{\b x_{(4)}^1 }  p(\b x_{\ell,i}^{(4)}), \quad \ell=1,\dots, 6, \ i=1,\dots, 123, }
where $\partial{\b x_{(4)}^1}$ denotes the order one derivatives in the normal direction,
   and   $\b x_{\ell,i}^{(4)}$ are internal face Lagrange  
      points described above.
      
On a second layer of triangle of this tetrahedron,  
  we have only 2D $P_{18}$ degrees of freedom left, at most,
  ranging from Bernstein index $(1,2,3,3,8,16)$ to index $(1,2,3,3,16,8)$ on the first edge of the
   triangle.
Thus, on this triangle, we have 2D 
\an{\label{2d-2} R^{(2)}_{18,-4} = \dim P_{18}-3 \dim P_4=145} degrees of freedom as
  the index $(1,2,3,3,7,17)$, $(1,2,3,3,6,18)$,  $(1,2,3,3,5,19)$, $(1,2,3,3,4,20)$ and $(1,2,3,3,3,21)$ are
   covered by the vertex degrees of freedom. 
That is,  we need only the following  
  degrees of freedom for the normal derivatives on the first triangle of the first tetrahedron
      of this 4D-simplex,
\a{ \partial{\b x_{(4)}^1 }  p(\b x_{\ell,i}^{(4)}), \quad \ell=1,\dots, 6, \ i=124,\dots, 268, }
where $\partial{\b x_{(4)}^1}$ denotes the order one derivatives in the normal direction,
   and   $\b x_{\ell,i}^{(4)}$ are internal face Lagrange  
      points described above.

Starting from the third layer to the sixth level of triangles, 
  we have one less vertex chopping from $(1,2,4,2,*,*)$ to $(1,2,4,5,*,*)$ so that we have sequentially
\an{\label{2d-3} \ad{ &\quad \ R^{(2)}_{20,-5} + R^{(2)}_{19,-4} + R^{(2)}_{18,-3} + R^{(2)}_{17,-2} \\
     & = 168+165+160+153 } } degrees of freedom.
The last few levels do not have vertex chopping but have 2 lines of edge chopping,
   $(1,2,*,2,2,*)$ and $(1,2,*,2,3,*)$ (here a skip of index implies the restriction is from
     edges of other tetrahedra), so that we can sum them together as
\an{\label{2d-4} \ad{ &\quad \ \sum_{i=8}^{16} R^{(2)}_{i,-1}  
       = C_{16+3,3}-C_{7+3,3}-(16-7)\cdot 3\cdot 3 = 768 } } degrees of freedom.
Together, by \eqref{2d-1}--\eqref{2d-4}, we have the number of degrees of freedom on the first tetrahedron
   (of index $(1,2,2,*,*,*)$--$(1,2,16,*,*,*)$), 
\an{\label{s4-1-1} \t{dof}_{\t{4s},1,1}  = 1682. }

On the second tetrahedron  (of index $(1,3,2,*,*,*)$--$(1,3,16,*,*,*)$), repeating above work, we skip the
  first layer of triangles and add one more layer of triangles at the end, i.e., we have 
\an{\label{s4-1-2} \ad{ \t{dof}_{\t{4s},1,2} &=R^{(2)}_{18,-4}+\sum_{i=1}^5 R^{(2)}_{15+i,-i} 
       +  \sum_{i=7}^{15} R^{(2)}_{i,-0}\\
  &=145+168+165+160+153+144\\
   &\qquad +(C_{15+3,3}-C_{6+3,3}-(15-6)\cdot 3)\\
    &  = 1640. } }

From the third tetrahedron to the sixth tetrahedron, we have a transitional pattern that we have specially  
\an{\label{s4-1-3} \ad{ \t{dof}_{\t{4s},1,3 } &= \sum_{i=2}^5 R^{(2)}_{15+i,-i} + \sum_{i=0}^1 R^{(2)}_{15+i,-i} 
       +  \sum_{i=6}^{14} C_{i+2,2}, \\
          \dots & \\
          \t{dof}_{\t{4s},1,6} &= \sum_{i=2}^2 R^{(2)}_{15+i,-i} + \sum_{i=0}^1 R^{(2)}_{15+i,-i} 
       +  \sum_{i=2}^{14} C_{i+2,2}.     } }  
We have a perfect 4D-simplex index set on the rest tetrahedron that 
\an{\label{s4-1-7} \ad{ \t{dof}_{\t{4s},1,7 } &=  \sum_{i=0}^1 R^{(2)}_{15+i,-i} 
       +  \sum_{i=2}^{14} C_{i+2,2}, \\
        \t{dof}_{\t{4s},1,8 } &=  \sum_{i=0}^0 R^{(2)}_{15+i,-i} 
       +  \sum_{i=2}^{14} C_{i+2,2}, \\
        \t{dof}_{\t{4s},1,9 } &= \qquad\qquad\qquad \ 
          \sum_{i=2}^{14} C_{i+2,2}, \\
          \dots & \\
          \t{dof}_{\t{4s},1,15} &= \qquad\qquad\qquad \ \sum_{i=2}^{8} C_{i+2,2}, \\
    \sum_{i=7}^{15} \t{dof}_{\t{4s},1,i} &= C_{20,4}-C_{11,4}-(16-7)\cdot 4\cdot C_{4,3} = 4371. } } 
The number of degrees of freedom for the normal derivative inside the 4D-simplex is, combining 
   \eqref{s4-1-1}--\eqref{s4-1-7},
\an{\label{4s1} \ad{ \t{dof}_{\t{4s},1} &=1682+1640+ 1547+ 1400+1250+1100+4371 \\
                      &= 12990. }  }

By \eqref{s4-1-1} and \eqref{4s1}, the number of degrees of freedom on all 4D-simplex is
 \an{\label{s4} \t{dof}_{\t{4s}} = 6\cdot (6965+12990) = 119730. }

 \section{The 5D-simplex degrees of freedom }

 As the function is $C^1$-continuous on the 5 face-4D-simplex of a 5D-simplex,
  supposedly we need to post $\dim P_{33-6\times 2}=\dim P_{21}$ 
  on the internal 5D Lagrange nodal degrees of freedom.
   
On the first triangle of this tetrahedron of this 4D simplex of the 5D-simplex,
  we have an additional edge of index $(2,2,2,1,*,*)$ being covered by edge degrees of freedom from
    triangles from other tetrahedra.
(Again, a reversed order index indicates an overlapping from other side.)
That is, we have a 2D $P_{21-3}=P_{18}$ polynomial on the triangle, instead of a 2D $P_{21}$ function.
Additionally the 5 Bernstein index from $(2,2,2,3,7,17)$ to $(2,2,2,3,3,21)$ are
       covered by the degrees of freedom at a vertex, i.e., chopping off a 3D $P_4$ corner. 
That is, on the first triangle, we have the following degrees of freedom,
\a{  \t{dof}_{\t{5s},0,0,1} = R^{(2)}_{18,-4} = C_{18+2,2} - 3 C_{4+2,2} =145. }
On the next triangle, we have only 1 edge-index $(2,2,3,1,*,*)$ instead of 3 edge-index overlapping.
Thus, we have 2 higher polynomial degrees and a $P_5$ corner chopping that  
\a{  \t{dof}_{\t{5s},0,0,2} = R^{(2)}_{20,-5} = C_{20+2,2} - 3 C_{5+2,2} =168. }
On the next triangle, we have no edge-index overlapping from its own tetrahedron, but
   1 edge-index $(  2,  2,  4,  8, 16,  1,)$ overlapping from other tetrahedron that
\a{  \t{dof}_{\t{5s},0,0,3} = R^{(2)}_{19,-5} = C_{19+2,2} - 3 C_{4+2,2} =165. }
On the next triangles, we have 
\a{  \t{dof}_{\t{5s},0,0,4} &= R^{(2)}_{18,-3}=160, \\
     \t{dof}_{\t{5s},0,0,5} &= R^{(2)}_{17,-2}=153,\\
     \t{dof}_{\t{5s},0,0,6} &= R^{(2)}_{16,-1}=144,\\ 
     \t{dof}_{\t{5s},0,0,7} &= R^{(2)}_{15,-0},\\ 
     \t{dof}_{\t{5s},0,0,8} &= R^{(2)}_{14,-0}=C_{14+2,2}-3 C_{2,2},\\
     \t{dof}_{\t{5s},0,0,9} &= R^{(2)}_{13,-0}=C_{13+2,2}-3 C_{2,2},\\
        \dots &  \\ 
     \t{dof}_{\t{5s},0,0,15} &=R^{(2)}_{7,-0}=C_{7+2,2}-3 C_{2,2},\\
     \sum_{i=7}^{15}\t{dof}_{\t{5s},0,0,15} &= C_{15+3,3}-C_{6+3,3}-3(15-6)=705. }
Together, on the first tetrahedron,  we have
\a{  \t{dof}_{\t{5s},0,1} &= 145+168+165+160+153+144+705=1640. }

On the second tetrahedron, we have
\a{ \t{dof}_{\t{5s},0,2} &=\sum_{i=0}^5 R^{(2)}_{15+i,-i} + \sum_{i=6}^{14} C_{i+2,2} \\
           &= 168+165+160+153+144+133+C_{14+3,3}-C_{5+3,3}\\
            & =1547. }
Which can be combined as 
\a{ \t{dof}_{\t{5s},0,2} &= R^{(3)}_{20,-5}=C_{20+3,3}-4 C_{5+3,3} =1547. }

Adding the first level and the rest levels of tetrahedral degrees of freedom, we get the
  degrees of freedom on the first 4D-simplex, 
\a{  \t{dof}_{\t{5s},1}&= 1640+\sum_{i=1}^5 R^{(3)}_{15+i,-i}+\sum_{i=7}^{15} R^{(3)}_{i,-0} \\
       &= 1640+ C_{20+4,4}-5 C_{5+4,4}-C_{6+3,3}-(15-6)(4)\\
       &  = 11520.  }
Repeating the count on the second 4D-simplex (first 4D-simplex has a different pattern),
     we get 
\a{  \t{dof}_{\t{5s},2}&= \sum_{i=1}^5 R^{(3)}_{15+i,-i}+\sum_{i=6}^{14} C_{i+3,3} \\
       &= R^{(4)}_{20,-5}=  C_{20+4,4}-5 C_{5+4,4}  = 9996.  }
Thus, we can combine the counting of degrees of freedom on all rest 4D-simplexes as a 5D-simplex counting
  that 
\an{\label{dof5} \ad{ \t{dof}_{\t{5s}}&= 11520 + \sum_{i=1}^5 R^{(4)}_{15+i,-i}+\sum_{i=7}^{14} R^{(4)}_{i,-0} \\
       &= 11520 +  R^{(5)}_{20,-5}-R^{(4)}_{6,-0} -(14-5)(5) \\
        & = 11520 +  C_{20+5,5}-6 C_{5+5,5}-C_{4+6,4}-(14-6)(5)\\
       &= 62888.  } }
We conclude the construction of $C^1$-$P_{33}$ in 5D by checking the total degrees of freedom 
   and the finite element space dimension, by  \eqref{1-s}, \eqref{2-s}, 
        \eqref{tri},  \eqref{tet}, \eqref{s4} and \eqref{dof5},
 \a{ \t{dof}_{\t{all}} &= \t{dof}_{\t{vert}} + \t{dof}_{\t{edge}}+ \t{dof}_{\t{tri}}  
          +  \t{dof}_{\t{tet}} + \t{dof}_{\t{4s}}    +   \t{dof}_{\t{5s}}
    \\ &=  122094+ 47520 + 57820 +  91890+ 119730+ 62888  \\
       &= 501942 = C_{33+5,5} = \dim P_{33} \  \t{in 5D}. }

\end{document}